\documentclass[12pt]{amsart}
\usepackage{amsfonts}
\usepackage{xcolor}
\usepackage{soul}
\usepackage{lscape}
\usepackage{enumerate}

\newtheorem{lemma}[subsection]{Lemma}
\newtheorem{theorem}[subsection]{Theorem}

\newenvironment{theorem1}{\par\noindent{\textbf{Almost all orbits attain almost bounded values}}\,\,\em}{\rm}
\newenvironment{conjecture1}{\par\noindent{\textbf{Conjecture 1}}\,\,\em}{\rm}
\newenvironment{conjecture2}{\par\noindent{\textbf{Conjecture 2}}\,\,\em}{\rm}
\newtheorem{conjecture}[subsection]{Conjecture}
\newtheorem{question}[subsection]{Question}
\newtheorem{corollary}[subsection]{Corollary}
\newcommand{\Comment}[1]{{\small\rm
#1}{\color{red}}} 
\renewcommand{\Comment}[1]{}
\title{Push-forward of geometric distributions under Collatz iteration: Part 1}
\author{Mary Rees}

\def\ep{\noindent{\hfill $\Box$}}
\begin{document}
\numberwithin{equation}{subsection}

\maketitle

\begin{abstract} Two conjectures are presented. The first, Conjecture 1, is that the pushforward of a {\em{geometric distribution}} on the integers  under $n$ Collatz iterates, modulo $2^p$, is usefully close to  uniform distribution on the integers modulo $2^p$, if $p/n$ is small enough.  Conjecture 2 is that the density is bounded from zero for the incidence of both $0$ and $1$ for the coefficients in the dyadic expansions of  $-3^{-\ell }$ on all but an exponentially small set of paths of a {\em{geometrically distributed}} random walk on the two-dimensional array of these coefficients. It is shown that Conjecture 2 implies Conjecture 1. At present, Conjecture 2 is unresolved.\end{abstract}

\section{Introduction}

Like many of  the best mathematical problems, the origin of the Collatz Conjecture is hazy. I cannot do better than refer to the introductory articles of Jeffrey Lagarias in the volume of papers edited by him  \cite{Lag2}.  My colleague Terry (CTC) Wall, whose engagement  with this question has been important to me (if not completely voluntary), has confirmed that he and fellow mathematicians came across the  Conjecture  in Cambridge in the 1950's \cite{Lag2}. My personal obsession took hold during the covid pandemic. In April 2020,  I was looking online for advice to pass on, and I came across Terence Tao's paper \cite{Tao}. This paper is a very important motivation for much of what follows, and so will be referred to very frequently, but first, here are a couple of observations, facile but worth note. 

The Collatz Conjecture is about an iteration on the integers, a dynamical system for which the phase space is the set of integers. These are fascinating but generally impossible exercises. Collatz' notebooks from the 1930's onwards contain many interesting examples (\cite{Lag1}, \cite{Lag2}), most of them still unsolved. Dynamicists always want to work with a topological space, preferably compact, or locally compact. A common technique, when studying Collatz iteration, is to use the embedding of the integers into the ring of dyadic integers  -- or triadic integers. Triadic integers  appear, largely implicitly, in \cite{Tao}. Dyadic integers appear in this paper, for essentially the same reason as triadic in \cite{Tao}. Dyadic and triadic  will appear together in Part 2. It is possible to extend Collatz iteration to the ring of dyadic integers \cite{Lag1}, and this embedding will be used here. But we will not be using the larger dynamical system. That is basically much simpler -- just because it is a standard shift space, a typical chaotic dynamical system.

Another tool which has become ubiquitous across mathematics is probability. Probability can be added in any context, and one can ask how ``most'' mathematical objects behave. If a problem seems intractable, put a probability on the phase space and try to prove the conjecture with high probability, preferably with probability approaching one. There is a long history of introducing probability into the Collatz Conjecture, back to the 1970's.  Stochastic models were constructed by Lagarias and Weiss \cite{Lag-Weiss}. \cite{Tao} is a more recent and extremely powerful application. There might be some hope that probabilistic arguments  not only give interesting new results but, ultimately, resolve the original conjecture. Tao passes on an argument from Ben  Green that a ``probability one'' result for the Collatz Conjecture could be as difficult to prove as the original conjecture. The current work grew out of an attempt to do something for the integers that get ``left behind'' by Tao's method. At the time of writing this remains a string of conjectures,  each implled by  the next,  but there does seem to be potential for further progress. The first and final conjectures in the sequence will be given shortly , as Conjecture 1 and Conjecture 2, appearing later as Corollary \ref{4.2} to  Conjecture \ref{4.1}, and  Conjecture \ref{4.13}. 
\subsection{Definitions}\label{1.1}

$Col:\mathbb Z\to \mathbb Z$ is defined by 
\begin{equation}\label{1.1.1}Col (N)=\begin{array}{ll}3N+1&\mbox{ if }N\mbox{ is odd,}\\
\frac{1}{2}n&\mbox{if }N\mbox{ is even}\end{array}\end{equation}
So $Col(0)=0$ and $Col $ preserves the sets of strictly positive integers $\mathbb N+1$ and strictly negative integers.
$Syr:\mathbb \mathbb Z\setminus \{ 0\} \to 2\mathbb Z+1$ is defined by
\begin{equation}\label{1.1.2}Syr(N)=N_2=Col^{x+y+1}(N)\end{equation}
if  $N=2^yN_1$ with $N_1$ odd and $Col(N_1)=2^{x}N_2$ and $N_2$ is odd. Then $Syr (1)=1$ and $Syr(-1)=-1$. It is convenient to extend $Syr$ to all nonzero integers, not just the odd ones, although of course $\mbox{Im}(Syr)=2\mathbb Z+1$, and $Syr $, too, preserves the sets of strictly positive, and strictly negative, integers. If $a_i$ ($i\ge 1$) is a sequence of numbers then we write (following \cite{Tao})
$$a_{[1,n]}=\sum _{i=1}^na_i.$$
If $N\in 2\mathbb Z+1$ then the sequence $x_i=x_i(N)$ of strictly positive integers  associated to $N$ is  defined by
$$Syr^i(N)=Col ^{x_i+1}(Syr^{i-1}(N))=Col^{x_{[1,i]}+i}(N).$$
\Comment{Another useful map to consider is 
$Col_2:\mathbb Z\to \mathbb Z$ is defined by 
\begin{equation}Col _2(n)=\begin{array}{ll}(3n+1)/2&\mbox{ if }n\mbox{ is odd,}\\
\frac{1}{2}n&\mbox{if }n\mbox{ is even}\end{array}\end{equation}
We then have 
$$Syr(n)=Col_2^{x_1}(n)$$
and
$$Syr ^i(n)=Col _2^{x_i}(Syr^{i-1}(n))=Col _2^{x_{[1,i]}}(n).$$}
If $n\ne 0$ is even and $n=2^{x_0}n_1$ with $n_1$ odd then for all $i\ge 1$,
$$Syr^i(n)=Col ^{x_{[0,i]}+i}(n)$$
\subsection{Collatz Conjecture}
There are several equivalent forms. The equivalence for the last form follows from the well-known  Theorem \ref{2.6}.
\begin{itemize} 
\item If $N\in \mathbb N+1$ then $Col ^k(N)\in \{1,2,4\} $ for all sufficiently large $k$.
\item if $n\in \mathbb N+2$ then for some $k$ we have $Col ^k(N)<N$.
\item If $N$ is a positive odd integer then $Syr^k(N)=1$ for all sufficiently large $N$.
\item If $N$ is a positive odd integer and $N>1$ then $Syr^k(N)<N$ for some $k$.
\item If $N$ is a positive odd integer with corresponding sequence $x_k$, then $x_k=2$ for all sufficiently large $k$.
\end{itemize}

For negative integers the conjecture is that all orbits of the Syracuse map $Syr$ end in one of three periodic orbits:
$$\begin{array}{l}-1\to -1;\ \ -5\to-7 \to -5;\\
  -17\to -25\to -37\to -55\to -41\to -61\to -91\to -17.\end{array}$$
In terms of the sequence $x_i$ the conjecture is that the sequence $x_i$ ends in the cycle $\dot{a}=aaa\cdots $ where $a$ is one of the following:
$$1:\ \ \ \ 12;\ \ \ 1112114.$$

The content of this Part 1 is that Conjecture 2 implies Conjecture 1. Definitions are given properly later.

\qquad\qquad

\begin{conjecture1}  Let $1<\mu<\infty $. Let $x_i$ be independent, identically distributed random variables on $\mathbb N+1$ with geometric distribution with mean $\mu $, that is, the probability that $x_i=m$ is $(\mu -1)^{m-1}\mu ^{-m}$. 

There is a constant $c>0$ which is bounded from $0$ is $\mu $ is bounded and bounded from $1$, such that the following holds for $k$ sufficiently large  and any $p\le ck$.   Let $N=N(x_1,\cdots x_k,1)$ be the smallest odd positive number with $x_i=x_i(N)$ for $i\le k$.

Then for any set $Y\subset \{ i\in\mathbb N:0\le i<2^p\} $, and for  $N=N(x_1,\cdots x_k,1)$, we have
$$|\mbox{Prob}(Syr^k(N)\in Y\mbox{ mod }2^p) |\le 4k\cdot 2^{-p}\#(Y).$$
\end{conjecture1}

\begin{conjecture2} The following holds for a constant $c>0$ depending on $\mu $ with  $1<\mu<\infty $. Fix $k$  and $p\le k$.  Let  the random variables $x_i$ be as in Conjecture 1. Let the dyadic expansion of $-3^{-j}$ be
$$-3^{-j}=\sum _{i=0}^\infty a_{i,j}2^i.$$
Consider the paths  $(\ell,x_{[\ell,k]}+p)$, $1\le \ell \le k$. Then apart from a set   of paths with probability $<e^{-ck}$, the density of $\ell\le k$ with $a_{x_{[\ell,k]}+p,\ell}=0$ is $\le (1-c)k$ and similarly for the density of the set of $\ell $ with $a_{x_{[\ell,k]}+p,\ell}=1$.\end{conjecture2}

\qquad\qquad

\section{Basic results}\label{2}

Versions of \ref{2.1}, \ref{2.2}, \ref{2.3} appear in publications  from the 1970's onwards (\cite{Crandall}, \cite{Kras}, \cite{Kras-Lag}). They are the basis of proofs of existence results for integers whose Collatz orbits terminate in the $1\mapsto 4\mapsto2\mapsto 1$ cycle or Syracuse orbits terminate in $1\mapsto 1$, in particular for results of the  type that the number of integers $N<N_0$ for which  the Syracuse  orbit  terminates in $1\mapsto 1$ is $\ge N_0^\alpha $ for some $\alpha >0$, for all sufficiently large $N_0$. These existence results are for integers $N$ such that the orbit  $Syr^k(N)$ is strictly decreasing in $k$, until $Syr^k(N)=1$. The difficulty with computing Syracuse or Collatz orbits is that they generally are not strictly decreasing. It is quite easy to see that for ``most '' $N$, in any reasonable sense, there will be $k$ for which $x_k=x_k(N)=1$ and in that case we have $Syr^k(N)>Syr^{k-1}(N)$. Generally if orbits increase, control is lost. Tao's very powerful methods \cite{Tao} cannot handle this. Nor indeed can the methods presented here, at least, not yet. 
\begin{lemma}\label{2.1} $$Syr (N)=\frac{3N+1}{2^k}\Leftrightarrow N=(-1)^kc_k+a_12^{k+1}$$
for some $a_1\in\mathbb N$ if $k$ is even, and $a_1\in \mathbb N+1$ if $k$ is odd, where
\begin{equation}c_k=\begin{array}{ll}2^{k-1}-2^{k-2}+\cdots +1=\sum _{j=0}^{k-1}(-1)^j2^j&\mbox{ if }k{\ is\ odd,}\\
\ &\ \\1+2^2+\cdots +2^{k-2}=\sum _{j=0}^{k/2-1}4^j&\mbox{ if }k\mbox{ is even.}\end{array}\end{equation}
Then
$$Syr(N)=6a_1+(-1)^k.$$
\end{lemma}
\noindent{\em{Proof}}
$$3c_k=2^k-(-1)^k$$
So
$$3(-1)^kc_k+1=(-1)^k2^k.$$
So if $N=(-1)^kc_k+a_12^{k+1}$ we have
$$\frac{3N+1}{2^k}=\frac{6a_12^k+(-1)^k2^k}{2^k}=6a_1+(-1)^k.$$
\ep

\begin{lemma}\label{2.2}
\begin{enumerate}[1.]
\item $Syr(4N-1)=6N-1$.
\item $Syr(8N+1)=6N+1$.
\item $Syr ^2(16N-5)=Syr(24N-7)=18N-5$.
\end{enumerate}
\end{lemma}
The proof of this is straightforward calculation.

\begin{lemma} \label{2.3}For all $k\ge 1$,
\begin{equation}\label{2.3.1}\mbox{Im}(Syr^k)=(6\mathbb N+1)\cup (6\mathbb N+5).\end{equation}
\end{lemma}
\noindent{\em{Proof}} It follows from Lemma \ref{2.1} (or \ref{2.2}) that (\ref{2.3.1}) holds for $k=1$. It will follow by induction if we can show that (\ref{2.3.1}) also holds for $k=2$. Using Lemma \ref{2.1}, it suffices to show that for any $k\in\mathbb N$, then for some even $j$
$$k2^{j+1}+c_j\equiv  \pm1\mbox{ mod }6$$
and for some odd $\ell $
$$k2^{\ell +1}-c_\ell \equiv \pm 1\mbox{ mod }6$$
For then by Lemma \ref{2.1} 
$$Syr(k2^{j+1}+c_j)=6k+1$$
and 
$$Syr(k2^{\ell +1}-c_\ell)=6k-1$$ 
But also by Lemma \ref{2.1}, $k2^{j+1}+c_j\in\mbox{Im}(Syr)$ and $k2^{\ell +1}-c_\ell \in\mbox{Im}(Syr)$.
First we consider the case of even $j=2t$.  We have
$$k2^{2t+1}+c_{2t}=k2^{2t+1}+1+2^2+\cdots +2^{2t-2}\equiv \pm 1\mbox{ mod }6$$
$$\Leftrightarrow k2^{2t}+2(1+\cdots 2^{2t-4})\equiv (\pm 1-1)/2\mbox{ mod }3$$
$$\Leftrightarrow k+2(t-1)\equiv (\pm 1-1)/2\mbox{mod }3$$
$$\Leftrightarrow 2k+t\equiv \pm 1\mbox{mod }3$$
Whatever $k$ is, there are infinitely many integers $t$ and $j=2t$ so that this holds.
Now we consider odd $\ell =2t-1$. Then  we need
$$k2^{2t}+1-2+2^2\cdots +2^{2t-2}\equiv  \pm1\mbox{ mod }6$$
$$\Leftrightarrow k2^{2t}+1(1+2\cdots +2^{2t-3})\equiv (\pm 1-1)/2\mbox{ mod }3$$
$$\Leftrightarrow 2k+2+t-2\equiv \pm 1\mbox{ mod }3$$
$$\Leftrightarrow 2k+t\equiv \pm 1\mbox{mod }3$$
which is the same equation as before. Once again, given $k$, there are infinitely many $t$ and $\ell =2t-1$ such that this holds. 

\ep

Most of us will have, at some stage, computed by hand  Collatz orbits for numbers up to (say) $2^{10}$, and will have used results such as \ref{2.4} and \ref{2.5}.

\begin{lemma}\label{2.4}For any $i$ and $j$, if $x_{[1,i-1]}(p)< j$ and $x_{[1,i]}(p)\ge j$ and $0<p< 2^j$  is odd, then \begin{equation}\label{2.4.1}\begin{array}{l}
x_\ell (p+2^j)=x_\ell (p)\mbox{  for }\ell \le i-1\\
x_{[1,i]}(p+2^j)\ge j\end{array}\end{equation}
 Furthermore: 
 \begin{equation}\label{2.4.2}\begin{array}{l}
 x_{[1,i]}(p+2^j)=j\mbox{ if }x_{[1,i]}(p)>j\\
 x_{[1,i]}(p+2^j)>j\mbox{  if }x_{[1,i]}(p)=j\end{array}\end{equation}
 In particular $x_{[1,i]}(p+2^j\ge j+1$  if $x_{[1,i]}(p)=j$.\end{lemma}
\noindent{\em{Proof}} $x_i(p)$ is the least integer $\ell \ge 1$ such that 
$$Col ^{1+\ell} (Syr ^{i-1}(p))=Col ^{1+x_{[1,i-1]}+\ell }(p)$$ is odd. By induction we have, for $1\le t\le i-1$
$$Syr^t(p+2^j)=Col ^{x_t+1}(Syr^{t-1}(p+3^{t-1}2^{j-x_{[1,t-1]}})$$
$$=Col^{x_t}(3Syr^{t-1}(p)+1+3^{t}2^{j-x_{[1,t-1]}})$$
$$=Syr^t(p)+3^t2^{j-x_{[1,t]}}$$
and hence)
\begin{equation}\label{2.4.3}Syr ^t(p+2^j)=Syr ^t(p)+3^t2^{j-x_{[1,t]}}\end{equation}
Similarly 
\begin{equation}\label{2.4.4}Col ^{1+j-x_{[1,i-1]}}(Syr^{i-1}(p+2^j)=Col^{1+j-x_{[1,i-1]}}(Syr^{i-1}(p))+3^i\end{equation}
where the right-hand side is odd if $x_{[1,i]}>j$, and even if $x_{[1,i]}=j$.
It follows that
\begin{equation}\label{2.4.5}
x_{[1,i]}(p+2^j)\begin{array}{lll}=j&\mbox{ if }&x_{[1,i]}(p)>j\\
>j&\mbox{ if }&x_{[1,i]}(p)=j\end{array}\end{equation}
\ep

\begin{lemma}\label{2.5}
\begin{enumerate}[1.]
\item $2^m-1$ and $2^{m-1}-1$ have the same Syracuse orbit if $m$ is even.
\item $2^m-1$ and $2^m +2^{m-1}-1$ have the same Syracuse orbit for all $m\ge 2$
\item $2^m +2^{m-3}-5$ and $2^m-5$ have the same Syracuse orbit for all $m\ge 6$.
\end{enumerate}
\end{lemma}
\noindent{\em{Proof.}}
There are many variations on these statements. The first is certainly well known. Unfortunately there are nowhere near enough common orbits for this to be the basis for an induction.

1. $$Syr^m(2^m-1)=(3^m-1)/2^q$$
where $2^{q-2}$ is the highest power of $2$ which divides $m$ --- so that $2^q$ is the highest power of $2$ dividing $3^m-1$ ---  and
$$Syr ^{m-1}(2^{m-1}-1)=(3^{m-1}-1)/2,\ \ Syr((3^{m-1}-1)/2)=(3^m-1)/2^q,$$
where $q$ is as before.

2. $$Syr ^{m-2}(2^m+2^{m-1}-1)=2^2\cdot 3^{m-2}+2\cdot 3^{m-2}-1,$$
and so
$$\begin{array}{l}Syr ^{m-1}(2^m+2^{m-1}-1)=(2^2\cdot 3^{m-1}+2\cdot 3^{m-1}-2)/2^{q+1}\\
\ \\
=(3^m-1)/2^q=Syr^m(2^m-1),\end{array}$$
where $2^q$ is the highest power of $2$ dividing $3^m-1$.

3. For each $p\ge 1$ with $3p< n$, 
$$Syr^{2p}(2^n-5)=3^{2p}\cdot 2^{n-3p}-5,$$
and 
$$\begin{array}{l}Syr ^{2p-2}(2^n+2^{n-3}-5)=Syr^{2p-2}(3^2\cdot 2^{n-3}-5)\\
\ \\
=3^{2+2p-2}\cdot 2^{n-3-3(p-1)}-5
=3^{2p}\cdot 2^{n-3p}-5=Syr ^{2p}(2^n-5).\end{array}$$
\ep

The following well-known and  very important identity is proved by straightforward induction.

\begin{theorem} \label{2.6}
The following are equivalent.
\begin{enumerate}[1.]
\item 
$N=2^{x_0}N_1$ for  $N_1\in \mathbb N$ odd, and $Syr^i(N)=Col^{x_{[0,i]}+i}(N)$ for $1\le i<k$ and $Syr ^k(N)=Col^{r+1}(Syr^{k-1}(N))$ for $r\ge x_k$.
\item
\begin{equation}\label{2.6.1}N\equiv -\sum _{i=0}^k3^{-i-1}2^{x_{[0,i]}}\mbox{ mod }2^{x_{[0,k]}}.\end{equation}
\end{enumerate}
Then 
\begin{equation}\label{2.6.2}Col ^{x_{[0,k]}+k}(N)\equiv \frac{3^k}{2^{x_{[0,k]}}}N+\sum _{j=1}^k\frac{3^{k-j}}{2^{x_{[j,k]}}}\mbox{ mod }3^k\end{equation}
\end{theorem}

\subsection{The numbers $N(x_1,\cdots x_k)$}\label{2.7}
\
\par

We will usually consider (\ref{2.6.1}) for $N$ odd that is $x_0=0$. We will  write $N=N(x_1,\cdots x_k)$ for the unique odd integer $N<2^{x_{[1,k]}}$ such that (\ref{2.6.1}) holds. This then identifies the odd positive integers $<2^n$ with finite sequences $(x_1,\cdots x_k)$  (for varying $k$) such that $x_{[1,k]}=n$.  Given $N=N(x_1,\cdots x_k)$ one can, of course, extend to an infinite sequence $x_i(N)$, where $x_i=x_i(N)$ for $i\le k-1$ and $x_k(N)\ge x_k$ and $x_i(N)$ is determined by $N$ for all $i$, that is, by $(x_1,\cdots x_k)$. We also note that $N=N(x_1,\cdots x_k,1)$ is the unique odd integer $N$ with $0<N<2^{x_{[1,k]}+1}$ and $x_i(N)=x_i$ for all $1\le i\le k$ The Collatz conjecture is that $x_i=2$ for all sufficiently large $i$.  This identification of an odd integer with a sequence $x_i$ of strictly positive integers is then the basis of all probabilistic results and conjectures. It is natural to regard the $x_i$ as independent random variables, usually identically distributed. Like \cite{Tao} we will focus on {\em{geometric distributions}}.
, which will be defined in the next section.

According to Lagarias \cite{Lag1} the first integer $k$ such that $3^k2^{-x_{[1,k]}}<1$ is known as the {\em{coefficient stopping time}} and apparently it is unknown whether the coefficient stopping time coincides with the {\em{stopping time}}, that is, the first $k$ for which $Syr^k(N)<N$.  Such an estimate can be  used for  the famous probabilitistic result, first known to be proved independently by Terras and Everett (\cite{Terras}, \cite{Everett}) that, for a positive density set of numbers $N<2^n$, approaching $1$ as $n\to \infty $, the stopping time is $<n$. A slightly more sophisticated estimate was used by Allouche \cite{Allouche} to prove that for a density approaching $1$, for various $\theta $, and for some least $k$, we have $Syr^k(N)<N^\theta $. The difficulty  is always in proving that this happens again. (But this was done in \cite{Tao} to a remarkable extent.) The following, for example, can be used to get it once. 

\begin{lemma}  \label{2.8} Let $0<\alpha <1/2$. Let $N=N(x_1,\cdots x_k,1)$. Suppose that 
\begin{equation}\label{2.8.1}\begin{array}{l}N\ge 2^{(1-\alpha )x_{[1,k]}}\\
3^{-k}2^{x_{[1,k]}}\ge 2^{\alpha x_{[1,k]}}\\
3^{-j}2^{\alpha x_{[1,k]}}\le 2^{x_{[j,k]}}\mbox{ for }1\le j\le k.\end{array}\end{equation}
Then
\begin{equation}\label{2.8.2}Syr^k(N)\le 2^{-\alpha x_{[1,k]}}(1+k)N<N^{1-\alpha /2} \end{equation}
assuming $k$ is large enough given $\alpha $.
 \end{lemma}
 \noindent{\em{Proof.}} We have, from (\ref{2.6.2}) 
 \begin{equation}\label{2.8.3}\begin{array}{l}Syr^k(N)=3^k\cdot 2^{-x_{[1,k]}}\left( N+\sum _{j=1}^{k}2^{x_{[1,j-1]}}\cdot 3^{-j}\right) \\

 \le 3^k\cdot 2^{-x_{[1,k]}}N\left( 1+\sum _{j=1}^{k}3^{-j}2^{\alpha x_{[1,k]} -x_{[j,k]}}\right) \\
 \le 2^{-\alpha x_{[1,k]}}N(1+k)
 \end{array}\end{equation}
 \ep

\subsection{Dyadic integers}\label{2.9}
The ring of dyadic integers is the completion of the ring of integers using the metric $d_2(n,m)=|n-m|_2=2^{-r}$, where $2^r$ is the highest power of $2$ dividing $n-m$. Every dyadic integer has a unique expansion 
$$\sum _{i=0}^\infty a_i2^i$$
with $a_i\in\{ 0,1\} $ for all $i$. The natural numbers $\mathbb N$ are precisely the dyadic integers which have a finite expansion of this form. All integers and all rational numbers with odd denominator are dyadic integers. For example
$$-1=(1-2)^{-1}-=\sum _{i=0}^\infty 2^i$$ and
$$-3^{-1}=(1-2^2)^{-1}=\sum_{i=0}^\infty 2^{2i}$$
All odd integers are units in the ring of dyadic integers, and $2$ is the only prime, modulo units. Rather trivially, the ring of dyadic integers is a unique factorisation domain. So the highest power of $2$ dividing a dyadic integer has meaning. Therefore the formulae for the Collatz map and the Syracuse map in (\ref{1.1.1}) and (\ref{1.1.2}) extend to the ring of dyadic integers. The ring is also the inverse limit of the rings $\mathbb Z/2^p\mathbb Z$. The Collatz and Syracuse maps  do not  extend to the quotient domain $\mathbb Z/2^p\mathbb Z$ --- unless, as \ref{2.4} shows, we quotient the range by $3^k\mathbb Z$ (for $Syr ^k$, for example). Nevertheless, if we keep the domain as the integers,  we can quotient the image by $2^p\mathbb Z$. We will use the following formula $Col ^{x_{[1,k]}+k}(N)$ as a dyadic integers, for $N=N(x_1,\cdots x_k)$ as in \ref{2.7}. We can then project down to $\mathbb Z/2^p\mathbb Z$ for any $p$. 

\begin{theorem}\label{2.10} Let $N=N(x_1,\cdots x_k)$ be an odd integer as in \ref{2.7}. As a dyadic integer
\begin{equation}\label{2.10.1}Col ^{x_{[1,k]}+k}(N)=Col ^{x_k+1}(Syr^{k-1}(N))=-3^k\left( \sum _{i=1}^{k}\sum _{m=0}^\infty a_{m+x_{[i,k]},i}2^m+A\right) \end{equation}
where $A=A(N)$ is an integer with $0\le A\le k$ and, as a dyadic integer,
$$-3^{-i}=\sum _{j=0}^\infty a_{j,i}2^j.$$
Now let $N=N(x_1,\cdots x_k,1)$. Then as a dyadic integer
\begin{equation}\label{2.10.2}
Syr^k(N)=-1+3^{-1}-3^k\left(\sum _{i=1}^k\sum _{m=0}^\infty a_{m+x_{[i,k]},i}2^m+A_1\right)\end{equation}
where $A_1$ is an integer with $-k\le A_1\le 2k+2$. 
\end{theorem}
\noindent{\em{Proof}}
By (\ref{2.6.1}), remembering that $x_0=0$ as $N$ is odd,
\begin{equation}\label{2.10.3}N\equiv -\sum _{i=0}^{k-1}3^{-i-1}2^{x_{[1,i]}}\mbox{ mod }2^{x_{[1,k]}},\end{equation}
where $x_{[1,0]}=0$ and in general $x_{[i,j]}=0$ if $i>j$. This implies that
\begin{equation}\label{2.10.4}N=\sum _{i=0}^{k-1}\{ -3^{-i-1}\} _{x_{[i+1,k]}}2^{x_{[1,i]}}-A2^{x_{[1,k]}}\end{equation}
where, if $y$ is an integer $\mbox{ mod }2^\ell$ for any $\ell \ge m$, then  $\{ y\}_m$ is the integer with $0\le \{ y\} _m<2^m$ and $y\equiv \{ y\}_m\mbox{ mod }2^m$. Then $A=A(N)$ is an integer and since  $0<N<2^{x_{[1,k]}}$ we have $0\le A\le k$. From (\ref{2.10.4}) we have
\begin{equation}\label{2.10.5}Col ^{x_1+1}(N)=\sum _{i=1}^{k-1}3\{ -3^{-i-1}\} _{x_{[i+1,k]}}2^{x_{[2,i]}}+(1+3\{- 3^{-1}\} _{x_{[1,k]}})2^{-x_1}-3A2^{x_{[2,k]}}\end{equation}

By induction we have, for $j<k$
\begin{equation}\label{2.10.6}Col ^{x_{[1,j]}+j}(N)=\begin{array}{l}\sum _{i=j}^{k-1}3^j\{ -3^{-i-1}\} _{x_{[i+1,k]}}2^{x_{[j+1,i]}}\\
\ \\
+\sum _{i=1}^{j}3^{j-i}(1+3^i\{ -3^{-i}\} _{x_{[i,k]}})2^{-x_{[i,j]}}
-3^jA2^{x_{[j+1,k]}}\end{array}\end{equation}
 and finally, for $j=k$, we have
\begin{equation}\label{2.10.7}Col ^{x_{[1,k]}+k}(N)=\sum _{i=1}^{k}3^{k-i}(1+3^i\{ -3^{-i}\} _{x_{[i,k]}})2^{-x_{[i,k]}}-3^kA\end{equation}
Now $3^{-i}$ is a dyadic integer. 
If we write
\begin{equation}\label{2.10.8}-3^{-i}=\sum _{m=0}^\infty a_{m,i}2^m\end{equation}
for $a_{m,i}\in \{ 0,1\} $ for all $m$, then we always have $a_{0,i}=1$ and 
\begin{equation}\label{2.10.9}\{ -3^{-i}\} _{x_{[i,k]}}=\sum _{m=0}^{x_{[i,k]}-1}a_{m,i}2^m.\end{equation}
So as a dyadic integer
\begin{equation}\label{2.10.10}3^{k-i}(1+3^i\{ -3^{-i}\} _{x_{[i,k]}})=-3^k(-3^{-i}-\{ -3^{-i}\} _{x_{[i,k]}})=-3^k\sum _{m=x_{[i,k]}}^\infty a_{m,i}2^{m-x_{[i,k]}}.\end{equation}
So interpreting all quantities as dyadic integers, (\ref{2.10.7}) becomes
 (\ref{2.10.1})
 
 Now we apply this to $N=N(x_1,\cdots x_k,1)$, so that $x_{k+1}=1$. Then we have
 \begin{equation}\label{2.10.11}\begin{array}{l}
 Col ^2(Syr^k(N))=-3^{k+1}\left( \sum _{i=1}^{k+1}\sum _{m=0}^\infty a_{m+x_{[i,k]}+1,i}2^m+A\right) \\
 \ \\
 =\displaystyle{\frac{3}{2}}(Syr^k(N)+1).\end{array}\end{equation}
 where now $0\le A\le k+1$. Now multiplying through by $2/3$, we obtain
 \begin{equation}\begin{array}{l}\label{2.10.12}Syr^k(N)=-1-3^k\left( \sum _{i=1}^k\sum _{m=1}^\infty a_{m+x_{[i,k]},i}2^m+\sum _{m=1}^\infty a_{m,k+1}2^m+2A\right) \\
 \ \\
 =-1+3^{-1}-3^k\left(  \sum _{i=1}^k\sum _{m=0}^\infty a_{m+x_{[i,k]},i}2^m-\sum _{i=1}^ka_{x_{[i,k]},i}+2A\right) 
 \ \\
 =-1+3^{-1}-3^k\left(  \sum _{i=1}^k\sum _{m=0}^\infty a_{m+x_{[i,k]},i}2^m+A_1\right) \end{array}\end{equation}
 for an integer $A_1$ with $-k\le A_1\le 2(k+1)$. This gives (\ref{2.10.2}), as required.
 \ep
 \qquad\qquad
 
The quantity $A=A(N)$ in  (\ref{2.10.1}) could be determined more precisely, if not exactly. As they stand, (\ref{2.10.1}) and (\ref{2.10.2}) are not precise formulae.. But they will be useful. The now-classical probabilistic results of Terras 
and Everett show that for some $\beta <1$ (which will be made more precise in Section \ref{3}), for all odd $1<N<2^n$ outside a set with at most $2^{\beta n}$ elements, we have $Syr^k(N)<N$ for some $k<n$. If we could use induction on this result we could make progress on the Collatz Conjecture. But we do not know how $Syr^k$ transforms sets of integers of large positive density. Tao \cite{Tao} makes substantial advance on this question. The statements of his results are quite involved, for very good reasons. 

 We note that if $k\le Cp$ and  and $X\subset \{ x:0\le x<2^p\} $ is a set with $\le 2^{\beta  p}$ elements for some $\mu<1$ then $\{ x+A\mbox{ mod }2^p:x\in X,0\le A\le k\} $ has $\le (k+1)2^{\beta  p}\le (Cp+1)2^{\beta  p}=o(2^p)$ elements. So although  the imprecise formula (\ref{2.10.2}) cannot show that Syracuse orbits are approximately uniformly distributed mod $2^p$ (for example) it might be useful for showing the image of a sparse set is sparse mod $2^p$. 

\section{Geometric distributions}\label{3}
\subsection{$x_i$ sequences with fixed mean}\label{3.1}

Write $N=N(x_1,\cdots x_k,1)$. \ref{2.8} gives conditions, including $3^k2^{-x_{[1,k]}}<2^{_\alpha x_{[1,k]}}$ under which 
$$Syr^j(N)<N^{1-\alpha /2}.$$

 If $x_{[1,k]}\le k\log 3/(\log 2(1-\alpha ))$ then this does not happen. Such sequences are not extremely rare. Now we consider  sequences $x_i$  for which $x_{[1,k]}=n$ for fixed $k$ and $n$.  We know from the generating function
$$x^k(1-x)^{-k}$$ 
that there are $(-1)^k\displaystyle{\left(\begin{matrix}-k\\n-k\end{matrix}\right)}=\displaystyle{\left(\begin{matrix}n-1\\k-1\end{matrix}\right)}$ such sequences $(x_1,\cdots x_k)$ and hence the same number of integers $N(x_1,\cdots x_k)$. Now write
$$\mu =\frac{n}{k}$$
Then 
$$1\le \mu \le n$$
We know that the most common value of $\mu $ -- fixing $n$ and varying $k$ --- is $2$.

Tao \cite{Tao} defined {\em{geometric distributions}}. For $0<\lambda <1$, the probability distribution $\mathbb G_\lambda $ on $\mathbb N+1$ (not Tao's notation) is defined by
$$\mathbb G_\lambda (X=m)=\lambda ^{m-1}(1-\lambda ),\ m\in \mathbb N+1.$$
This probability distribution has mean
$$(1-\lambda )\sum _{m=1}^\infty m\lambda ^{m-1}=(1-\lambda )\frac{d}{d\lambda }\left( \sum _{m=1}^\infty \lambda ^m\right) $$
$$=(1-\lambda )\frac{d}{d\lambda }\frac{\lambda }{1-\lambda }=(1-\lambda )\left(\frac{1}{1-\lambda}+\frac{\lambda }{(1-\lambda )^2}\right) =\frac{1}{1-\lambda }$$
So if the mean is $\mu $ we have
$$\lambda =1-\mu ^{-1}$$
and the geometric probability distribution with mean $\mu $ is $\mathbb G_{1-\mu ^{-1}}$. This probability distribution extends to a joint distribution on $k$ independent random variables, for each $k$, and hence to a probability measure on $\cup _{k=1}^\infty (\mathbb N+1)^k $ which we will also denote by $\mathbb G_{1-\mu ^{-1}}$. We denote by $[x_1,\cdots x_k]$ the set, or {\em{cylinder}} defined by 
$$\bigcup _{m\ge k}\{ [y_1,\cdots y_m]:y_i=x_i,1\le i\le k\} .$$
Then the probability measure $\mathbb G_{1-\mu ^{-1}}$ is defined by 

 $$\mathbb G_{1-\mu ^{-1}}([x_1\cdots x_k])=\prod _{i=1}^k\mu ^{-1}(1-\mu ^{-1})^{x_i-1}$$
 $$=(\mu-1)^{-k}\left( \frac{\mu -1}{\mu }\right) ^{x_{[1,k]}}=\mu ^{-k}(1-\mu ^{-1})^{x^{[1,k]}-k}$$
 So all $k$-cylinders $[x_1,\cdots x_k]$ with $x_{[1,k]}=n$, for fixed  $k$ and $n$, have the same measure. If $\mu =2$, then all $k$-cylinders, for varying $k$ and fixed $n$, have the same measure. If we fix $k$, 
 $\mathbb G_{1-\mu ^{-1}}$ can also be regarded as a a probability measure on the set of odd integers $N(x_1,\cdots x_k,1)$, by identifying $N(x_1,\cdots x_k,1)$ with the cylinder $[x_1,\cdots x_k]$. For fixed $n=x_{[1,k]}$, for varying $k$, this is a set of $2^{n-1}$ odd positive integers $N<2^{n+1}$, that is, exactly half of the odd positive  integers $<2^{n+1}$. We note that all the odd positive integers $<2^{n+1}$ are of the form $N(x_1,\cdots x_{k+1})$ for some $k\ge 0$ and $x_{[1,k+1]}=n+1$. If $\mu=2$ the mass assigned by $\mathbb G_{1-\mu ^{-1}}$ to each integer $N(x_1,\cdots x_k,1)$ with $x_{[1,k]}=n$ is $2^{-n}$. Without a bound on $x_{[1,k]}$, the set of integers $N(x_1,\cdots x_k,1)$ is infinite because, for each $i$, $x_i$ can take any positive integer value. However,  there are advantages in considering this infinite set, because if we fix $\mu $ then the $\mathbb G_{1-\mu ^{-1}}$ -measure of the cyllinders $[x_1,\cdots x_k,1]$ is concentrated on those $[x_1,\cdots x_k,1]$ with $|x_{[1,k]}-\mu k|\le C\sqrt{k}$, for $C$ sufficiently large. 
 
 Of course one can also identify the cylinder $[x_1,\cdots x_k]$ with the coset of integers which are equal to $N(x_1,\cdots x_k,1)\mbox{ mod }2^{x_{[1,k]}+1}$. That is more the approach of \cite{Tao}, but is more appropriate when restricting to  $\mu >\log 3/\log 2$.
  
  We recall the following completely standard lemma, which is no surprise. 
 
 \begin{lemma}\label{3.2} If $\mu >1$ and $k/n\sim \nu <1$ then 
 \begin{equation}\label{3.2.1}\begin{array}{l}
 \left(\begin{matrix}n-1\\k-1\end{matrix}\right)\mu ^{-k}(1-\mu ^{-1})^{n-k}\\
 \ \\
 \sim\displaystyle{\frac{1}{\sqrt{2\pi n\nu (1-\nu )}}}\exp n(\nu \log(\mu ^{-1}\nu ^{-1})+(1-\nu )\log(1-\mu ^{-1})(1-\nu)^{-1}).\end{array}\end{equation}
 Moreover the maximum of this is attained at $\nu =\mu ^{-1}$, when the exponential term on the right becomes $1$. 
 \end{lemma}
 \noindent {\em{Proof.}} We repeat the standard proof, using Stirling's formula
 $$m!\sim \sqrt{2\pi }m^{m+1/2}e^{-m}.$$
 So 
 $$\left(\begin{matrix}n-1\\k-1\end{matrix}\right)\mu ^{-k}(1-\mu ^{-1})^{n-k}$$
 $$\sim \frac{\sqrt{n-1}}{\sqrt{2\pi(k-1)(n-k)}}\exp n\left(\frac{k}{n}\left( \log\frac{n}{k}+\log \mu ^{-1}\right)+\frac{n-k}{n}\left(\log \frac{n}{n-k}+\log (1-\mu ^{-1})\right) \right) $$
 $$\sim \frac{1}{\sqrt{2\pi n\nu (1-\nu)}}\exp n(\nu \log (\nu ^{-1}\mu ^{-1})+(1-\nu )\log( (1-\nu )^{-1}(1-\mu ^{-1}))).$$
 To compute the maximum in $\nu $, write
 $$g(\nu )=g_\mu (\nu )=\nu \log\frac{\mu ^{-1}}{\nu }+(1-\nu )\log\frac{1-\mu ^{-1}}{1-\nu }.$$
 Then
 $$g'(\nu)=\log \frac{\mu ^{-1}}{\nu }-\log\frac{1-\mu ^{-1}}{\nu}-1+1=0\ \Leftrightarrow\ \frac{\mu ^{-1}}{\nu}=\frac{1-\mu ^{-1}}{1-\nu}$$
 $$\Leftrightarrow \nu =\mu ^{-1}$$
 We have
 $$g''(\nu )=-\frac{1}{\nu}-\frac{1}{1-\nu }<0$$
 and so $\nu =\mu ^{-1}$ gives a maximum. 
 \ep
 
 Applying \ref{3.2} with $\mu =2$ we have the following.
 \begin{corollary}\label{3.3} If $k/n\sim \nu $ then
 \begin{equation}\label{3.3.1}\begin{array}{l}\left(\begin{matrix}n-1\\k-1\end{matrix}\right)\sim 2^n\displaystyle{\frac{1}{\sqrt{2\pi n\nu (1-\nu )}}}\exp n(\nu \log(\nu ^{-1}/2)+(1-\nu )\log((1-\nu)^{-1}/2)\\
 \ \\
 = 2^n\displaystyle{\frac{1}{\sqrt{2\pi n\nu (1-\nu )}}}2^{-n\lambda (\nu)}\end{array}\end{equation}
 with $\lambda (\nu )\le0$ and bounded from $0$ if $0<\nu <1$ and $\nu $ is bounded from $1/2$.\end{corollary}
 
 From \ref{3.3} the number of cylinders $[x_1,\cdots x_k]$ with $x_{[1,k]}=n$ and $3^k<2^{n(1-\varepsilon)}$ for a sufficiently small $\varepsilon $ is $<2^ne^{-\lambda n}$ for a $\lambda >0$. By \ref{2.6} this is also the number of integers $N(x_1,\cdots x_k)$ with $x_{[1,k]}=n$ --- that is, $0<N(x_1,\cdots x_k)<2^{n}$ --- and $3^k<2^{n(1-\varepsilon)}$. Hence by \ref{2.8}, as $n\to \infty $, the proportion of $N<2^n$ for which $Col ^{x_{[1,k]}+k}(N)<N$ for some $k<n$  is at least $1-ce^{-\lambda n}$ for  constants $\lambda >0$ $c>0$.
  
 Lemma \ref{3.2} also gives a direct proof of the Central Limit Theorem for the probability measure/joint distribution $\mathbb G_{1-\mu ^{-1}}$. More accurate information can be obtained by using the Fourier transform method applied to this example. The Central Limit Theorem shows (of course) that the measure $G_{1-\mu ^{-1}}$ is concentrated on those cylinders $x_{1,k]}$ --- equivalently on those integers $N(x_1,\cdots x_k,1)$ --- such that $|x_{[1,k]}-k\mu|\le C\sqrt{k}$, where,  as $C\to \infty$, the measure outside this set tends to $0$ for all sufficiently large $k$ given $C$.
 \begin{theorem}\label{3.4} Let $A\subset [-L,L]\subset \mathbb R$ be an interval. Then for a constant $C$ depending only on $L$, 
 $$\left\vert \mathbb G_{1-\mu ^{-1}}\left( \frac{x_{[1,k]}-k\mu }{\sqrt{k}}\in A\right) -\frac{1}{\sqrt{2\pi \mu (\mu-1)}}\int _Ae^{-t^2/(2\mu(\mu -1))}dt\right\vert \le Ck^{-1/50}$$
 \end{theorem}
 The proof is redacted.
 
\Comment{ \noindent{\em{Proof}}
 This is a  standard proof of the Central Limit Theorem applied to our example. We identify the Fourier transform of the probability distribution of $(x_{[1,k]}-k\mu )/\sqrt{k}$ to some degree of accuracy. 
  The probability distribution of $x_{[1,k]}-k\mu $ is the $k$-fold involution $*^k\mathbb G_{\mu, 1-\mu ^{-1}}$ where $\mathbb G_{\mu , 1-\mu ^{-1}}$ is the probability distribution of $x_1-\mu$, that is, 
 $$\mathbb G_{\mu, 1-\mu ^{-1}}({m-\mu })=\mu ^{-1}(1-\mu ^{-1})^{m-1}.$$
 Now
 $${\widehat{ \mathbb G}}_{\mu ,1-\mu ^{-1}}(t)=\sum _{m=1}^\infty e^{i(\mu -m)t}\frac{1}{\mu }\left( 1-\frac{1}{\mu}\right) ^{m-1}$$
 $$=\frac{e^{i\mu t}}{\mu -1}\sum _{m=1}^\infty e^{-imt}\left(\frac{\mu -1}{\mu }\right) ^m$$
 $$=\frac{e^{i(\mu -1)t}}{\mu }\left( 1-\frac{\mu -1}{\mu}e^{-it}\right) ^{-1}.$$
 So if $\mathbb G=*^k\mathbb G_{\mu, 1-\mu ^{-1}}$, we have
 $${\widehat {\mathbb G}}(t)=({\widehat{ \mathbb G}}_{\mu ,1-\mu ^{-1}}(t))^k$$
 $$=\frac{e^{ikt(\mu -1)}}{\mu ^k}\left( 1-\frac{\mu -1}{\mu}e^{-it}\right) ^{-k}$$
 $$=\left( \frac{e^{i(\mu -1)t}}{\mu -(\mu -1)e^{-it}}\right) ^k=\left( \frac{e^{i\mu t}}{\mu e^{it}-\mu +1}\right) ^k.$$
 Now 
 $$|\mu e^{it}-\mu +1|^2=\mu ^2\sin^2t+(\mu \cos t-(\mu -1))^2=\mu ^2+(\mu-1)^2-2\mu (\mu-1)\cos t$$
 $$=\mu ^2+(\mu-1)^2-2\mu (\mu-1)+2\mu (\mu-1)(1-\cos t)=1+2\mu (\mu -1)(1-\cos t)\ge 1$$
 with equality if and only if $\cos t=1$. If $\cos t\le 1-\delta $ then we see that 
 $$|{\widehat {\mathbb G}}(t)|\le (1+2\mu (\mu-1)\delta )^{-k}$$
 If $t=o(1)$ then
 $${\widehat {\mathbb G}}(t)=\left( \frac{1+i\mu t-\mu ^2t^2/2-i\mu ^3t^3/6\cdots }{1+i\mu t-\mu t^2/2-i\mu t^3/6\cdots }\right) ^k$$
 $$=\left( \frac{1-\mu ^2t^2/2(1+i\mu t)^{-1}-i\mu ^3t^3/6(1+i\mu t)^{-1}\cdots }{1-\mu t^2/2(1+i\mu t)^{-1}-i\mu t^3/6(1+i\mu t)^{-1}\cdots }\right) ^k$$
 $$=\left( \frac{1-\mu ^2t^2/2+i\mu ^3t^3/2-i\mu ^3t^3/6\cdots }{1-\mu t^2/2-i\mu ^2t^3/2-i\mu t^3/6\cdots }\right) ^k$$
 $$=(1-\mu (\mu -1)t^2/2+i\mu t^3\cdots )^k=\exp -k(\mu (\mu -1)t^2/2+i\mu t^3\cdots )$$
 If $t-2m\pi =o(1)$ then 
 $$e^{ik\mu t}=e^{2\pi ikm\mu }e^{k\mu(t-2\pi im )}$$
 and 
  $${\widehat {\mathbb G}}(t)=e^{2\pi ikm\mu }\exp -k(\mu (\mu -1)t^2/2+i\mu t^3\cdots )$$
  
 Now the Fourier transform of the probability distribution for $x_{[1,k]}-k\mu )/\sqrt{k}$ is $\widehat{\mathbb G}(t/\sqrt{k})$, and we now have, if $t/\sqrt{k}=o(1)$, then 
$$\widehat{\mathbb G}(t/\sqrt{k})=\exp (-\mu (\mu -1)t^2/2+i\mu t^3/\sqrt{k}\cdots ).$$
In particular, this is valid for $|t|\le k^{1/10}$.

Now we need to use the inverse Fourier transform. We will use what feels like the most basic definition, for any probability measure $\mathbb G$, for any bounded measurable function $g$ with support in a finite interval $[-L,L]$
$$\int g(x)d\mathbb G(x)=\lim _{\Delta \to\infty }\frac{1}{2\pi}\int _{-L}^{L } \int _{-\Delta }^{\Delta }e^{itu}g(u){\widehat {\mathbb G}}(t)dt\ du$$
$$=\lim _{\Delta \to\infty }\frac{1}{2\pi}\int _{-L }^{L}\int _{-\Delta }^{\Delta }e^{itu}g(u)\int _{-\infty }^{\infty }e^{-itx}d\mathbb G(x)\ dt\ du$$
$$=\frac{1}{2\pi }\int _{-\infty }^{\infty }\int _{-L }^L \lim _{\Delta \to \infty }g(u)\int _{-\Delta }^{\Delta } e^{-it(u-x)}dt\ du\ d\mathbb G(x)$$
$$=\frac{1}{\pi }\int _{-\infty }^{\infty }\int _{-L}^L \lim _{\Delta \to \infty }g(u)\frac{\sin \Delta (u-x)}{u-x}du\ d\mathbb G(x)$$
The change of order of integration is justified because  two of the integration limit sets are finite and $\mathbb G$ is a probability measure. 
Meanwhile, by standard methods of contour integration,
$$\frac{1}{\pi} \int _{-L}^L\frac{\sin\Delta (u-x)}{u-x}du=\frac{1}{\pi }\int _{x-\Delta L}^{x+\Delta L}\frac{\sin v}{v}dv=1+O(\Delta ^{-1/2}$$
So now we consider
$$\frac{1}{\pi } \int _{-L}^L\frac{g(u)-g(x)}{u-x}\sin \Delta (u-x)du$$
Our aim is to show that this  is bounded by $O(\Delta ^{-\alpha }$ for some suitable $\alpha $. 
Because our integral involves the Dirichlet kernel $\sin (\Delta y)/y$, we will initally consider $g$ being piecewise $C^2$. We will then sandwich the characteristic function $\chi A$ , for our interval $A$,  between two piecewise $C^1$ functions, one of which is $0$ outside $A$ and $1$ on a amaller  interval with endpoints within $\varepsilon $ of those of $A$, and the other of which is $1$ on $A$ and $0$ outside a slightly larger interval. again with endpoints within $\varepsilon $ of those of $A$. We will assume as we may do that support of both our functions $g$ is contained in $[-L,L]$. First we note that if $h$ is piecewise $C^1$ with support in $[-L,L]$ and $h'$ denotes the first derivative of $h$ and $\Vert . \Vert$ denotes uniform norm then
$$\left\vert \int h(v)\sin \Delta v)dv\right\vert \le 4L\Vert h\Vert (\delta \Delta ) ^{-1}+2L\delta \Vert h'\Vert $$
This is obtained by splitting up the support of $h$ into intervals of length $\delta $ and approximating $h$ to within $\delta \Vert h'\Vert $ on each of these intervals by a constant function. 
Now we want to apply this with 
$$h(u)=\begin{array}{lll}\frac{g(u)-g(x)}{u-x}& if &u
\ne x\\
\ &\ &\ \\
g'(x)& if &u
= x\end{array}$$
If $g$ is piecewise $C^2$ then $h$ is piecewise $C^1$ with 
$$h'(u)=\begin{array}{lll}\frac{g'(u)(u-x)-(g(u)-g(x))}{(u-x)^2}& if &u
\ne x\\
\ &\ &\ \\
g'(x)& if &u
= x\end{array}$$
We then have
$$\Vert h\Vert \le \Vert g\Vert +\Vert g'\Vert ,$$
$$\Vert h'\Vert \le \Vert g'\Vert +\Vert g''\Vert .$$
So now we have
$$\left \vert \frac{1}{\pi } \int _{-L}^L\frac{g(u)-g(x)}{u-x}\sin \Delta (u-x)du\right\vert \le 4L(\Vert g\Vert +\Vert g'\Vert) (\delta \Delta ) ^{-1}+2L\delta (\Vert g'\Vert +\Vert g''\Vert )$$
We now suppose that $g$ is one of our functions sandwiching $\chi _A$, so that $\Vert g\Vert \le 1$, $\Vert g'\Vert =O(\varepsilon ^{-1})$ and $\Vert g''\Vert =O(\varepsilon ^{-2})$. So we then have 
$$\left \vert \frac{1}{\pi } \int _{-L}^L\frac{g(u)-g(x)}{u-x}\sin \Delta (u-x)du\right\vert \le C_1 (\delta \Delta \varepsilon) ^{-1}+\delta \varepsilon ^{-2}).$$
We then take $\delta =\varepsilon ^3$ and $\Delta\ge \varepsilon ^{-5}$. This gives 
$$\left \vert \frac{1}{\pi } \int _{-L}^L\frac{g(u)-g(x)}{u-x}\sin \Delta (u-x)du\right\vert \le C_2\varepsilon $$
for both choices of $g$. So for both choices of $g$  and for any $\Delta \ge \varepsilon ^{-5}$ we have 
$$\left\vert \int g(x)d\mathbb G(x)-\frac{1}{2\pi}\int _{-L}^{L } \int _{-\Delta }^{\Delta }e^{itu}g(u){\widehat {\mathbb G}}(t)dt\ du\right \vert \le C_3\varepsilon .$$
Moreover this estimate works with $\mathbb G$ being  {\em{any}} probability measure because we have a uniform bound on the integrand. 
However we have seen that if $\Delta \le k^{-1/10}$ then
$$\left\vert \frac{1}{2\pi}\int _{-L}^{L } \int _{-\Delta }^{\Delta }e^{itu}g(u)({\widehat {\mathbb G}}(t)-\exp (-\mu (\mu -1)t^2/2)dt\ du\right\vert \le C_4\Vert g\Vert k^{-2/5}.$$
We deduce that if $\mathbb Q$ is the normal distribution with density function 
$$1/\sqrt{2\pi \mu (\mu -1)}\exp (-t^2/(2\mu (\mu -1))),$$
then
$$\left\vert \int g\ d\mathbb G-\int g\ d\mathbb Q\right\vert \le 2C_3\varepsilon + C_4\Vert g\Vert k^{-2/5}$$
$$\le C_5k^{-1/50}$$
taking $\Delta =\varepsilon ^{-5}$
Since $\chi _A$ is sandwiched between our two choices of $g$ and the difference between the integrals of the two choices of $g$ with respect to $\mathbb Q$ is $O(\varepsilon)$, we obtain our result. 
\ep
}

The map $N\mapsto Syr^k(N)$ is many-to-one in general and of course that is the whole point of the Collatz conjecture. However, it is possible that $N\mapsto Syr^k(N)$ is injective restricted to the set of numbers for which the Syracuse stopping time $i>k$ where $i$ is the first integer (if it exists) with $Syr ^i(N)<N$.  \cite{Tao} gives a result of this type, which is at first sight surprising (Corollary  6.3), which holds for numbers which are generic for $\mathbb G_{1/2}$. We have the following for $\mathbb G_{1-\mu ^{-1}}$ for any $\mu <\log 3/\log 2$, which is an easy result of this type, rather far from what might be true. 

\begin{lemma}\label{3.5} For any $\mu <\log 3/\log 2$, let 
 $$M(\mu)=3^{-1}(1-2^\mu 3^{-1})^{-1}.$$
 Then for $N=N(x_1,\cdots x_k,1)$ and $N'=N(x_1',\cdots x_k',1)$ with $x_{[1,k]}=x'_{[1,k]}=n$ and  $x_{[1,i]}\le \mu i$ and $x'_{[1,i]}\le \mu i$ for all $i\le k$,  we have
 \begin{equation}\label{3.5.1}|Syr ^k(N)-Syr^k(N')|> 3^k2^{-n}(|N-N'|-M(\mu)).\end{equation}
 Consequently the   map $N\mapsto Syr ^k(N)$ is at most $\lfloor M(\mu )\rfloor $ -to-one  restricted restricted to the set of $N=N(x_1,\cdots x_k,1)$ with $x_{[1,k]}=n$ and  $x_{[1,i]}\le \mu i$ for all $i\le k$.\end{lemma}
 
\noindent{\em{Proof.}} We have 
$$Syr ^k(N)=\sum _{i=1}^k3^{k-i}2^{-x_{[i,k]}}+3^k2^{-x_{[1,k]}}N=3^k2^{-x_{[1,k]}}\left( N+\sum _{i=1}^k2^{x_{[1,i-1]}}3^{-i}\right) $$
So  since $x_{[1,k]}=x_{[1,k]}'=n$ we have
$$Syr ^k(N)-Syr ^k(N')=3^k2^{-n}(N-N'+\sum _{i=1}^{k}2^{x_{[1,i-1]}}3^{-i}-\sum _{i=1}^{k}2^{x_{[1,i-1]}'}3^{-i}). $$
Now 
$$\left\vert \sum _{i=1}^{k}3^{-i}(2^{x_{[1,i-1]}}-2^{x_{[1,i-1]}'})\right\vert <\sum _{i=1}^\infty 3^{-i}2^{(i-1)\mu }$$
$$=M(\mu ).$$
This gives (\ref{3.5.1}), and the final deduction is immediate.  The righthand-side of (\ref{3.5.1}) might of course be negative, but can only be so for at most $\lfloor M(\mu )\rfloor $ values of $N'$, given $N$. 
\ep

We remark that if $\mu =n_0/k_0$ for coprime  positive integers $n_0$ and $k_0$ then $M(\mu )$ can be replaced by $M'(\mu )$ where $M'(\mu )$ is defined as the maximum of all numbers of the form 
$$\sum _{i=2}^{k_0}2^{x_{[1,i-1]}}3^{-i}(1-2^{n_0}3^{-k_0})-\frac{2}{3}$$
where $x_{[1,n_0]}=k_0$, $x_{[1,i]}\le \mu i$ for all $i<n_0$. The $k$-tuple $(x_1,\cdots x_k)$ which achieves the maximum has the maximum value of $x_{[1,i-1]}$ for each $i\le k$. The subtraction of $2/3$ is because if $x_i'$ is the sequence for $N'$ then $x_i'\ge 1$ for all $i$ and
$$\sum _{i=1}^\infty 2^{i}3^{-i-1}=\frac{2}{9}\left( 1-\frac{2}{3}\right) ^{-1}=\frac{2}{3}.$$

\section{The Conjectures}\label{4}

We start with the main theorem of \cite{Tao}, Tao's Theorem 1.3.

\begin{theorem1}Let $f:\mathbb N+1\to \mathbb R$ be any function with $\lim _{N\to \infty }f(N)=+\infty $. Then one has $Col _{\mbox{min}}(N)<f(N)$ for almost all $N\in \mathbb N+1$.\end{theorem1}

\qquad\qquad

 For technical reasons ``almost all'' is in terms of what Tao calls ``logarithmic density'', rather than the normal density which weights all integers in a finite interval equally.  The proof uses a sophisticated induction, in order to show that, for a suitable  constant $\alpha >1$ for ``almost all '' $N$, there are successive iterates $Syr^{i_{[1,j]}}$ in the Syracuse orbit with $N_j=Syr^{i_{[1,j]}}(N)<N^{1/\alpha ^j}$, for $j\le J$, where $J$ can be taken arbitrarily large if $N$ is large enough. The first step  in the induction, with $j=1$ is Allouche's refinement \cite{Allouche} of the result of Terras and Everett . That removes  a proportion $N_0^{-\beta  }$ of numbers $<N_0$ permanently from consideration for a suitable $\beta >0$. After that first step it becomes much more difficult, because if is not clear how $Syr^{i_1}$ transforms sets of large positive density or very small density. These are numbers $\mbox{mod }3^n$ for some suitable $n$. Rather than looking at all numbers $\mbox{ mod }3^n$, Tao works with fibres of projections $\mbox{ mod }3^m$ for $m<n$: probably with good reason, although at first sight it seems a strange thing to do. It turns out that the Fourier transform of the key Proposition 1.14 can actually be solved. This then allows the inductive step to be completed. After the first step, in getting 
 $Syr ^{i_{j+1}}(N_j)<N_j^{1/\alpha }$, the proportion of logarithmic density lost is bigger, being $1/\log ^c(N_j)$ -- but that is enough, with the telescoping series. 
 
 It is an astonishing result. But as already mentioned in the introduction, at each stage numbers are lost and cannot be recovered. This seems to always be the way in dealing with Collatz iteration. From the second step of the induction, the numbers which are considered are those for which $x_{[1,k]}/k$ is close to $2$ At least for the very first step, following Terras, it is only the numbers with $x_{[1,k]}/k>\log 3/\log 2$ but after that it is close to $2$.

Our first conjecture is as follows. We use the notation of \ref{2.10}. In particular the $a_{j,i}$ are the coefficients in the dyadic expansion of $-3^{-i}$. 
\begin{conjecture}\label{4.1}
Define 
\begin{equation}\label{4.1.1}\begin{array}{l}
Syr_{k,p}(x_1,\cdots x_k)=-3^k\sum _{i=1}^k\sum _{j=0}^{p-1} a_{j+x_{[i,k]},i}2^j\\
\ \\
(Syr _{k,p})_*\mathbb G_{1-\mu ^{-1}}(X)=\mathbb G_{1-\mu ^{-1}}(\bigcup\{ [x_1,\cdots x_k]:Syr_{k,p}(x_1,\cdots x_k)\mbox{ mod }2^p\in X\}\\
\ \\
\mathbb S_{k,p,\mu}=(S_{k,p})_*\mathbb G_{1-\mu ^{-1}}\end{array}\end{equation}
That is,  $(Syr_{k,p})_*\mathbb G_{1-\mu ^{-1}}$ denotes the pushforward of the measure $\mathbb G_{1-\mu ^{-1}}$ under the map $(Syr)_{k,p}$
Then for any $1<\mu <\infty $ and  $0<c_1$ there exists $c_0>0$, depending on $c_1$ and $\mu $, bounded from $0$ if $\mu $ is bounded from $1$ and $c_1$ bounded from $0$ such that if $p\le c_0k$  and for any $0\le m<2^p$,
\begin{equation}\label{4.1.2}
|(\mathbb S_{k,p,\mu}(\{m\})-2^{-p}|\le c_12^{-p}\end{equation}
 \end{conjecture}

{\em{If Conjecture \ref{4.1} is true}}  then using Theorem \ref{2.10}  and Lemma \ref{2.8}, we  have  a result close to uniform distribution on the odd integers $\mbox{ mod }2^p$ of the iterates $Syr^k(N)$ --- not quite because the presence of the constant $A$ in \ref{2.10}.
\begin{corollary}\label{4.2}
Let $1<\mu <\infty$. Then if Conjecture \ref{4.1} is true,  for $c_0$ as in \ref{4.1} if $p\le c_0k$ then    for any set $Y\subset \{ i\in\mathbb N:0\le i<2^p\} $ we have  
\begin{equation}\label{4.2.1}\begin{array}{l}|\mathbb G_{1-\mu ^{-1}}([x_1,\cdots x_k]):Syr^k(N(x_1,\cdots x_k,1))\in Y)\\
\ \\
\le (1+c_1)\cdot (2k+2)2^{-p}\#(Y)\end{array}\end{equation}

Consequently, for a suitable universal constant $\alpha _1>0$,
outside a set of $\mathbb G_{1-\mu ^{-1}}$-measure $\le k2^{-\alpha _1p}$ we will have for $N=N(x_1,\cdots x_k,1)$ and for any $r$ with $p/4\le r\le p/3$
\begin{equation}\label{4.2.2}Syr^{k+r}(N)<2^{-r/5}Syr^k(N).\end{equation}\end{corollary}
\noindent{\em{Proof.}}
We apply \ref{2.8} with $N$ replaced by $N_1=Syr^k(N)\mbox{ mod }2^p$. More precisely, we apply (\ref{2.8.1}) to (\ref{2.8.3}) with $N$ replaced by $N_1$ and with $k$,  $x_{[1,k]}$ and $x_{[j+1,k]}$ replaced by $r$, $x_{[k+1,k+r]}$ and $x_{[k+j+1,k+r]}$, because $x_j(N_1)=x_{k+j}(N)$. In (\ref{2.8.2}) and (\ref{2.8.3}) we replace   $Syr^k(N)$ by $Syr^r(N_1)$. 
For  any fixed $p/4\le r\le p/3$ and $\alpha =1/5$ (for example) the following conditions are satisfied for $N_1$ outside a set of $\le 2^{p(1-\alpha _1)}$ elements, where throughout, this $x_j$ means $x_j(N_1)$, that is, $x_{k+j}(N)$ and $\alpha _1$ is a  constant depending only on $\alpha $.
\begin{equation}\label{4.2.3}\begin{array}{l}3r/2\le x_{[1,r]}\le 3r,\\
N_1\ge 2^{(1-\alpha )x_{[1,r]}},\\
3^{-r}2^{x_{[1,r]}}\ge 2^{\alpha x_{[1,r]}},\\
3^{-j}2^{\alpha x_{[1,r]}}\le 2^{[x_{[j,r]}}\mbox{ for }1\le j\le r.\end{array}\end{equation}
The last three conditions are simply the conditions of (\ref{2.8.3}) with the replacements indicated. We apply (\ref{3.2}) to see that for fixed $p\ge r\ge p/3$ and $\alpha =1/5$ these conditions hold for $N_1$  outside a set of $\le 2^{p(1-\alpha _1)}$ elements, where $\alpha _1$ is simply a universal constant given by \ref{3.2}.  For the penultimate condition of (\ref{4.2.3}) we note that $2^8>3^5$ and so $\log 3/\log 2<8/5$ and yet
$$3^{-r}2^{x_{[1,r]}}< 2^{x_{[1,r]}}/5\Rightarrow x_{[1,r]}<(5/4)\log 3/\log 2=2-(2-(5/4)\log 3/\log 2).$$
For the last condition of (\ref{4.2.3}) we can, for example, split into the cases $j\le r/2$ and $j\ge r/2$. If $j\le r/2$ then the condition is implied by 
$$2^{x_{[1,r]}/5}\le 2^{x_{[r/2,r]}}$$
and if $j\ge r/2$ we can use 
$$3^{-r/2}2^{x_{[1,r]}/5}\le 1.$$
\ep

All work to date has been on proving Conjecture \ref{4.1}. I believe progress has been made - producing a string of conjectures each implying the previous one, with a final conjecture which looks plausible and doable (to me) but which so far remains unresolved. Even if Conjecture 4.1 is proved, it is not, on its own, a useful analogue  of the Tao's method \cite{Tao} for his second inductive step for any value of $\mu $. This is because we have no lower bound on $c_0$ and  almost certainly have $p< <\mu k$ with current methods, if they can be made to work. By \ref{3.3}, we must have $p=O((\mu -1)\log (\mu -1)k)$.  We only have information (from \ref{4.2.1}) about $Syr ^k(N(x_1,\cdots x_k,1))\mbox{ mod }2^p$, when we start with $N(x_1,\cdots x_k,1)$ being any odd integer and if $\mu <\log 3/\log 2$ we will get 
$$Syr ^k(N)>N$$
 on a set of large $\mathbb G_{1-\mu ^{-1}}$-measure.
 
 It does  not look as if there is any particular advantage in having the constant $c_1$ in Conjecture \ref{4.1} being small, and one might as well choose $c_1=1$, for example. The question then is how small $c_0$ needs to be in terms of $\mu $ for fixed $c_1$.

\begin{question}\label{4.3}
How large can the constant $c_0$  of Conjecture \ref{4.1} be, depending on $\mu $, for fixed $c_1$? 
 \end{question}

\qquad
\qquad 
 
 Tao's proof of the main result in \cite{Tao} uses Fourier transform, which is natural to do in any probability proof. So now we consider the Fourier transform of $\mathbb S_{k,p,\mu}$. We have 
$${\widehat{\mathbb S }_{k,p,\mu}}(\xi )=\sum _{m=0}^{2^p-1}\mathbb S _{k,p,\mu}(\{m\} )e^{-2\pi im\xi /2^p}.$$
Then ${\widehat{\mathbb S }}_{k,p,\mu}(0)=1$, as is of course true for any probability measure. Also, $|{\widehat{\mathbb S}} _{k,p,\mu}(\xi )|\le 1$ for all $\xi $, and with strict inequality unless $\mathbb S _{k,p,\mu}$ is a point mass. Again, this fact about probability measure transforms is generally true. Now we have the following basic lemma.

\begin{lemma}\label{4.4} Conjecture \ref{4.1} is true if
\begin{equation}\label{4.4.1}\sum _{\xi \ne 0}|{\widehat{\mathbb S}} _{k,p,\mu}(\xi )|\le c_1.\end{equation}
\end{lemma}
\noindent{\em{Proof}} By Fourier inversion
$$\mathbb S _{k,p,\mu}(\{ m\} )=2^{-p}\sum _\xi {\widehat{\mathbb S}} _{k,p,\mu}(\xi )e^{2\pi im\xi /2^p}=2^{-p}+2^{-p}\sum _{\xi \ne 0}{\widehat{\mathbb S}} _{k,p,\mu}(\xi )e^{2\pi im\xi /2^p}.$$
\ep

So we give another conjecture, which, by \ref{4.4} would imply Conjecture \ref{4.1}. 
\begin{conjecture}\label{4.5} For some $c_2(\mu )>0$  for $k$ sufficiently large and $p\le k$, for all $0<\xi<2^p$
\begin{equation}\label{4.5.1}|{\widehat{\mathbb S}} _{k,p,\mu }(\xi )|\le 2^{-c_2k}.\end{equation}
\end{conjecture}
Conjecture \ref{4.5} looks a lot stronger than 1.17. of \cite{Tao},  just as \ref{4.1} looks a lot stronger than the result 1.14 of \cite{Tao}. Then if $c_0<c_2$ we have
$$\sum _{0<\xi<2^p}|{\widehat{\mathbb S}} _{k,p,\mu }(\xi )|\le 2^{(c_0-c_2)k}$$
Conjecture \ref{4.5} itself does not, therefore,  make much restriction on $p$. But to deduce Conjecture \ref{4.1} we need $p\le c_0k$ with $c_0<c_2$ and $k$ sufficiently large given $c_2-c_0$ and $c_1$. In the longer term it would  be better not to take modulus of ${\widehat{\mathbb S}} _{k,p,\mu }(\xi )$ but  for the moment the effort is concentrated on proving \ref{4.5}. If that is achieved, improvements will be sought. 

\subsection{}\label{4.6} A first thought is that there should be an inductive procedure for obtaining an estimate such as Conjecture \ref{4.5} on $|{\widehat{\mathbb S}} _{k,p,\mu}(\xi )|$. But the roughly  corresponding bound Proposition 1.17  of \cite{Tao} is  difficult and follows on from Tao's rather simple Lemma  6.2 which has no analogue in our case. However, we can consider what an inductive procedure would involve. We will fix  $\mu $ in what follows. We will also fix $p$ making only the restriction that $p\le k$.  We drop $\mu $ and $k$ as indices, so write $\mathbb S_{k}$ for $\mathbb S_{k,p,\mu }$ and $\mathbb G$ for $\mathbb G_{1-\mu ^{-1}}$. 
Write
\begin{equation}\label{4.6.1}m_i=m_i(x_{[i,k]})=\sum _{t=0}^{p-1}a_{t+x_{[i,k]},i}2^t\mbox{ mod }2^p\end{equation}
The Fourier transform ${\widehat{\mathbb S}}_{k}$, like the Fourier transform of any probability measure on $\mathbb Z/2^p\mathbb Z$, is a convex sum of $2^p$'th roots of unity. We have
\begin{equation}\label{4.6.2}{\widehat{\mathbb S}}_{k}(\xi )=\sum _{x_1,\cdots x_k}\mathbb G([x_1,\cdots x_k])\exp (-2\pi i\sum _{j=1}^km_j([x_{[j,k]}])\xi /2^p)\end{equation}

We note that $\mathbb G([x_1,\cdots x_k])$ depends only on $x_{[1,k]}$. Since $\mathbb G$ is a product measure this is a sum of products. We now want to write this as a sum of products of sums.
We can  break up into segments. Fix $r$ and $k_j$ such that $k_0=0$ and $k_j\ge 1$ for $1\le j\le r$ such that $k_{[1,r]}=k$. Now write  

 \begin{equation}\label{4.6.3}\begin{array}{l}{\widehat{\mathbb S }}_{k,p_{[j,r]},p_j}(\xi )=\\
 \ \\
\displaystyle{ \sum  \{ \prod _{\ell =k_{j-1}+1}^{k_j}\exp (-2\pi im_\ell(x_{[\ell,k]})\xi /2^p) \mathbb G([x_{k_{j-1}+1},\cdots x_{k_j}]):\  x_{[k_{j-1}+1,k_j]}=p_j\} }\\
 \ \\
= \displaystyle{\left(\frac{\mu -1}{\mu }\right)} ^{p_j}\mu ^{k_j-k_{j-1}}\sum \left \{\exp \left(-2\pi i \sum _{\ell =k_{j-1}+1}^{k_j}m_\ell (x_{[\ell,k ]})\xi/2^p  \right):\  x_{[k_{j-1}+1,k_j]}=p_j\right \} .\end{array}\end{equation}
Then we have
\begin{equation}\label{4.6.4}{\widehat{\mathbb S _{k}}}(\xi )=\sum _{p_1,\cdots p_r}\prod _{j=1}^r{\widehat{\mathbb S}} _{k,p_{[1,j-1]},p_j}(\xi ).\end{equation}
We also have 
\begin{equation}\label{4.6.5}\begin{array}{l}|{\widehat{\mathbb S }}_{k,p_{[j,r]},p_j}(\xi )|\le \\
\ \\
\sum\{  \mathbb G([x_{k_{j-1}}+1,\cdot ,x_{k_j}]):\ x_{[k_{j-1}+1,k_j]}=p_j\} =\left(\begin{matrix}p_j-1\\k_j-k_{j-1}-1\end{matrix}\right)  \left(\displaystyle{\frac{\mu -1}{\mu }}\right) ^{p_j}\mu ^{k_j-k_{j-1}}\end{array}\end{equation}

Write
\begin{equation}\label{4.6.6}\begin{array}{l}\lambda (\ell _1,\ell _2,q,u,\xi /2^p)=\\
\ \\
\left(\begin{matrix}q-1\\\ell _2-\ell _1-1\end{matrix}\right)^{-1} \sum \left \{\exp \left(-\displaystyle{\frac{2\pi i\xi}{2^p}} \sum _{\ell _1+1}^{\ell _2}m_\ell (x_{[\ell,\ell _2]+u})\right):\  x_{\ell _1+1}+\cdots +x_{\ell _2}=q\right \} 
\end{array}\end{equation}

Then we have , for fixed $k_1\cdots k_r$
\begin{equation}\label{4.6.7}\begin{array}{l}|{\widehat{\mathbb S _{k}}}(\xi )|\le \\
\ \\
\sum _{p_1,\cdots p_r}\left( \prod _{j=1}^r\lambda (k_{j-1},k_j,p_j,p_{[j+1,r]}, \xi/2^p )\cdot \left(\begin{matrix}p_j-1\\k_j-k_{j-1}-1\end{matrix}\right) \displaystyle{ \left(\frac{\mu -1}{\mu }\right) ^{p_j}}\mu ^{k_j-k_{j-1}}\right)\end{array}\end{equation}

Now we make the following conjecture. 
\begin{conjecture}\label{4.7} For $p\le k$ and  constants $c_3$, $c_4>0$ which are independent of $k$  and for fixed $k_1\cdots k_r$ with $k_j-k_{j-1}\le c_3^{-1}$ and any  $\xi \ne 0$, we have,   for at least $c_4r$  of the numbers $j\le r$:
\begin{equation}\label{4.7.1}|\lambda (k_{j-1},k_j,p_j,p_{[j+1,r]}, \xi /2^p)|\le e^{-c_4}\end{equation}
apart from a set $A$ of $(p_1,\cdots p_r)$ such that 
\begin{equation}\label{4.7.2}\sum _{(p_1,\cdots p_r)\in A}\prod _{j=1}^r\left(\begin{matrix}p_j-1\\k_j-k_{j-1}-1\end{matrix}\right)  \left(\frac{\mu -1}{\mu }\right) ^{p_j}\mu ^{k_j-k_{j-1}}\le e^{-c_4k}.\end{equation}
\end{conjecture}

Conjecture \ref{4.7} would then imply that 
\begin{equation}\label{4.7.3} {\widehat{\mathbb S}} _{k
}(\xi )\le e^{-c_3c_4^2k}+e^{-c_4k}\end{equation}
and hence implies Conjecture \ref{4.5} and hence Conjecture  \ref{4.1}. So this leads us to consider sums of the form
$$\sum _{j=1}^te^{2\pi i b_j}$$ 
where $0\le b_j<1$. We have the following frequently used lemma.

\begin{lemma}\label{4.8}
\begin{equation}\label{4.8.1}\left\vert \sum _{j=1}^te^{2\pi i b_j}\right\vert =\sqrt{t^2-4\sum _{j<\ell }\sin ^2\pi (b_j-b_\ell )}.\end{equation}
\end{lemma}
\noindent{\em{Proof.}} Simple calculation. 
\begin{equation}\label{4.8.2}\begin{array}{l}\left\vert \sum _{j=1}^te^{2\pi i b_j}\right\vert ^2=\sum _{j=1}^te^{2\pi i b_j}\overline{\sum _{\ell =1}^te^{2\pi i b_\ell }}\\
\ \\
=t+\sum _{j\ne \ell }e^{2\pi i(b_j-b_\ell )}\\
\ \\
=t+2\sum _{j<\ell}\cos 2\pi (b_j-b_\ell )\\
\ \\
=t^2-2\sum _{j<\ell }(1-\cos 2\pi (b_j-b_\ell))\\
\ \\
=t^2-4\sum _{j<\ell }\sin ^2\pi (b_j-b_\ell ).\end{array}\end{equation}
\ep
\subsection{}\label{4.9}
Now we apply \ref{4.8} to obtain a simpler condition to obtain (\ref{4.7.1})
{\em{for all $\xi \ne 0$}} in the case when $k_t-k_{t-1}$ and $p_t$ are bounded but $p_t-(k_t-k_{t-1})$ is sufficiently big. Since we expect $p_t-(k_t-k_{t-1})$ to be close to $(\mu -1)(k_t-k_{t-1})$ this means that we need $k_t-k_{t-1}$ to be sufficiently big, and we cannot expect better than $k_t-k_{t-1}=O((\mu -1)^{-1})$, that is $c_3^{-1}=O((\mu-1)^{-1})$ for $c_3$ as in Conjecture \ref{4.7}.  We assume that $\xi $ is odd (replacing $p$ by $p-u$ if $u$ is the maximum power of $2$ such that $2^u$ divides $\xi $) and by replacing $\xi $ by $-\xi $ if necessary, we assume that $0<\xi <2^{p-1}$. Now fix such a $\xi $.  If $m$ is any integer, write $\{ m\} $ for the integer $\mbox{ mod }2^p$ such that 
\begin{equation}\label{4.9.2}\begin{array}{l}1-2^{p-1}\le \{ m\} \le 2^{p-1},\\
m\equiv \{ m\} \mbox{ mod }2^p.\end{array}\end{equation}
Thus, we have
\begin{equation}\label{4.9.2}\begin{array}{l}|\{m\} |\le |m|,\\
|\{ m\} |\le 2^{p-1}\\
|\{ m_1+m_2\} |\le \{ m_1\} |+|\{ m_2\} |\\
|2^r\{ m\} |\le 2^r|\{ m\} |\end{array}\end{equation}
Fix $t$. Fix $x_{k_{t-1}},\cdots x_\ell,\cdots x_{k_t}$, where $x_{[\ell,k_t]}-(k_t-\ell )$ is sufficiently large   that $x_\ell $ can be varied by at least $n'$ for some suitable $n'$ (to be determined), keeping $x_i$ fixed for $k_{t-1}\le i\le k_t$ and $i\ne \ell $. Write $n=x_{[\ell,k]}$. Write
\begin{equation}\label{4.9.3}m(j+n,\ell)=\sum _{i=0}^{p-1}a_{i+n+j,\ell}2^i\end{equation}

From \ref{4.8}, $\lambda (k_{t-1},k_t,p_t,p_{[t+1,r]}, \xi /2^p)$ will be boundedly less than $1$ if we can show that for some fixed constant $C_0>0$, the following holds for some  $n'\le p_t-(k_t-k_{t-1})$, and some $\ell $ with  $k_{t-1}\le \ell <k_t$ and $x_{[k_{t-1}+1,\ell-1]}$ with 
$$x_{[k_{t-1}+1,\ell -1]}+n'+k_t-\ell\le p_t.$$
For some $j_1\ne j_2$  with $1\le j_1,j_2\le n'$
 \begin{equation}\label{4.9.4}|\{(m(j_1+n,\ell )-m(j_2+n,\ell))\xi \} |\ge C_02^p.\end{equation}
 
So now we find a simpler sufficient condition for (\ref{4.9.4})  
 We note that if $p$ is bounded then for $C_0<2^{-p}$,  if (\ref{4.9.4}) does {\em{not}} hold then $m(j_1+n,\ell)=m(j_2+n,\ell )$. The only way that can be true is if 
 \begin{equation}\label{4.9.5}\begin{array}{l}a_{n+i+j_1,\ell}=a_{n+i+j_2,\ell}\mbox{ for all }1\le j_1\le j_2\le n'\mbox{  and }0\le i<p,\mbox{ that is }\\
  a_{n+i+j,\ell}=a_{n+i+1,\ell}\mbox{ for all }1\le j\le n'\end{array}\end{equation}. So the case $p\le n'$ is reduced to a simpler condition. From now on we assume that $p$ is sufficiently large that (in particular) $p>n'$. We fix $\ell $, $\xi $, $n$.

Now for $1\le j\le n'$
\begin{equation}\label{4.9.6}\begin{array}{l}m(j+n,\ell)=\\
\ \\
\sum _{i=0}^{n'-1-j}a_{n+i+j,\ell}2^i+2^{n'-j}\sum _{i=0}^{p-1-n'}a_{n+n'+i,\ell}2^i+2^{p-j}\sum _{i=0}^{j-1}a_{n+p+i,\ell}2^i\\
\ \\
=b(j,n')+2^{n'-j}d+2^{p-j}c(j)\end{array}\end{equation}
So
\begin{equation}\label{4.9.7}b(n',n')=0,\ b(j,n')<2^{n'-j},\ c(j)<2^{j}.\end{equation}
 Then
\begin{equation}\label{4.9.8}m(j_+n,\ell)-m(n'+n,\ell))=b(j,n')+(2^{n'-j}-1)d+2^{p-n'}(2^{n'-j}c(j)-c(n'))\end{equation}
So
\begin{equation}\label{4.9.9}\begin{array}{l}(2^{n'-j_2}-1)(m(j_1+n,\ell)-m(n'+n,\ell))-(2^{n'-j_1}-1)(mj_2+n,\ell)-m(n'+n,\ell))\\
=f(j_1,j_2,n')+2^{p-n'}g(j_1,j_2,n')\end{array}\end{equation}
where
\begin{equation}\label{4.9.10}\begin{array}{l}f(j_1,j_2,n')=(2^{n'-j_2}-1)b(j_1,n')-(2^{n'-j_1}-1)b(j_2,n')\\
\ \\
g(j_1,j_2,n')=(2^{n'-j_2}-1)(2^{n'-j_1}c(j_1)-c(n'))-(2^{n'-j_1}-1)(2^{n'-j_2}c(j_2)-c(n'))\\
\ \\
=2^{2n'-j_1-j_2}(c(j_1)-c(j_2))+2^{n'-j_2}c(j_2)-2^{n'-j_1}c(j_1)+(2^{n'-j_1}-2^{n'-j_2})c(n').\end{array}\end{equation}

where both $f(j_1,j_2,n')$ and $g(j_1,j_2,n')$ are $O(2^{n'})$. In fact
\begin{equation}\label{4.9.11}|f(j_1,j_2,n')|<2^{2n'-j_1-j_2}\end{equation}
So if $|\{(m(j_1+n,\ell)-m (j_2+n,\ell))\xi \} |\le 2^{p-n''}$ for all choices of $j_1<j_2\le n'$, that is, (\ref{4.9.4}) does {\em{not}} hold, then it is also true that 
\begin{equation}\label{4.9.12}|\{ (f(j_1,j_2,n')+2^{p-n'}g(j_1,j_2,n'))\xi\} |\le 2^{p-n''+n'}\end{equation}
 for all choices of $j_1$ and $j_2$. We now assume that $n''$ is  bounded but somewhat larger than $n'$.

 Write 
 \begin{equation}\label{4.9.13}\begin{array}{l}\xi =\xi _1+2^{n'}\xi _2+2^{p-1-n''}\xi _4=\xi _3+2^{p-1-n''}\xi _4\mbox{ with}\\
 0<\xi _1<2^{n'},\ \ 0\le \xi _2<2^{p-1-n''-n'},\ \ 0\le \xi _4< 2^{n''},\ \ 0<\xi _3<2^{p-1-n''}\end{array}\end{equation}
 Then 
 \begin{equation}\label{4.9.14}\begin{array}{l}2^{p-n'}g(j_1,j_2)\xi \equiv 2^{p-n'}g(j_1,j_2)\xi _1\mbox{ mod }2^p\\ 
 \ \\
 |\{ f(j_1,j_2,n')\xi _3\} |<2^{p-1-n''+2n'-j_1-j_2}\end{array}\end{equation}
 So if (\ref{4.9.12}) holds we have
 \begin{equation}\label{4.9.15}|\{ f(j_1,j_2,n')\xi _42^{p-1-n''}+2^{p-n'}g(j_1,j_2,n'))\xi_1\} |<2^{p-n''+2n'-j_1-j_2}+2^{p-n''+n'}\end{equation}
that is
 \begin{equation}\label{4.9.16}|\{ f(j_1,j_2,n')\xi _4+2^{n''-n'+1}g(j_1,j_2,n'))\xi_1\}_{n''+1}|<2^{2n'-j_1-j_2+1}+2^{n'+1},\end{equation}
that is for some integer $A$ with 
$$|A|<2^{2n'-j_1-j_2+1}+2^{n'+1},$$
 we have
\begin{equation}\label{4.9.17}  f(j_1,j_2,n')\xi _4+2^{n''-n'}g(j_1,j_2,n')\xi_1\equiv A\mbox{ mod }2^{n''+1}\end{equation}
Now we consider this for choices of $j_1$, $j_2$, $n'$, $n''$

\subsection{$j_1=1$, $j_2=2$, $n'=3$, $n''=6$}\label{4.10}

We have 
\begin{equation}\label{4.10.1}\begin{array}{l}b(2,3)=a_{n+2,\ell},\ \ b(1,3)=a_{n+1,\ell}+2a_{n+2,\ell},\\
\ \\
  c(1)=a_{n+p,\ell},\\
\ \\
 c(2)=a_{n+p,\ell}+2a_{n+p+1,\ell}\\
\ \\
c(3)=a_{n+p,\ell}+2a_{n+p+1,\ell}+4a_{n+p+2,\ell}.\\
\ \\
f(1,2,3)=b(1,3))-3b(2,3)=a_{n+1,\ell}-a_{n+2,\ell}\\
\ \\
g(1,2,3)=8(c(1)-c(2))+2c(2)-4c(1)+2c(3)\\
\ \\
=-16a_{n+p+1,\ell}+2a_{n+p,\ell }+4a_{n+p+1,\ell}-4a_{n+p,\ell}+2a_{n+p,\ell}+4a_{n+p+1,\ell}+8a_{n+p+2,\ell}\\
\ \\
=-8a_{n+p+1,\ell}+8a_{n+p+2,\ell}\end{array}\end{equation}
So we have
\begin{equation}\label{4.10.2}f(1,2,3)\xi _4+2^3g(1,2,3)\xi _1\equiv (a_{n+1,\ell}-a_{n+2,\ell})\xi _4+2^6(a_{n+p+2,\ell}-a_{n+p+1,\ell})\xi _1\mbox{ mod }2^7\end{equation}

Then if (\ref{4.9.16}) holds and if 
 \begin{equation}\label{4.10.3}\begin{array}{l}N\equiv (a_{n+1,\ell}-a_{n+2,\ell})\xi _4+2^6(a_{n+p+2,\ell}-a_{n+p+1,\ell})\xi _1\mbox{ mod }2^7\\
 \ \\
0\le N<2^7\end{array}\end{equation}
 then since $|\{ N\} |_7<2^5$ that means either $0\le N<2^5$ or $2^5+2^6<N<2^7$. So the coefficients of $2^5$ and $2^6$ in the dyadic expansion of $N$ are both $0$ or both $1$.  Clearly one way in which this can happen is if $a_{n+p+1,\ell }=a_{n+p+2,\ell}$ and either  $0\le \xi _4<2^5$ or $a_{n+1,\ell }=a_{n+2,\ell }$.  Now suppose that $a_{n+p+1,\ell }\ne a_{n+p+2,\ell }$. Then since $-2^6\equiv +2^6\mbox{ mod }2^7$, we get a solution to (\ref{4.10.3}) only if the coefficient of $2^5$ in the dyadic expansion of $N_1$ is $1$ and the coefficient of $2^6$ is $0$ where $0\le N_1<2^6$ and $N_1\equiv (a_{n+1,\ell}-a_{n+2,\ell })\xi _4\mbox{ mod }2^7$. We need $a_{n+1,\ell}\ne a_{n+2,\ell }$ for this. But if $a_{n+1,\ell }=0$ and $a_{n+2,\ell }=1$ then
 \begin{equation}\label{4.10.4}\begin{array}{l}-2^6<(a_{n+1,\ell }-a_{n+2,\ell})\xi _4\le 0\\
 N_1=2^7-2^6<(a_{n+1,\ell }-a_{n+2,\ell})\xi _4\ge 2^6\end{array}\end{equation}
 So the only possibility is that $a_{n+1,\ell}=1$ and $a_{n+2,\ell }=0$. In addition, it is only possible if $2^5\le \xi _4<2^6$. 
 
 In summary, we have the following theorem.
 
 \begin{theorem}\label{4.11} Fix a suitable $M$ and 
 $$x_{[k_{t-1},k_t]}\ge M'\ge M^2+p_t-(k_t-k_{t-1}).$$
  Then there exists $0<\lambda _0=\lambda _0(M,M')<1$ such that the following holds for any odd integer $\xi $ with $0<\xi <2^{p-1}$.
Write  
\begin{equation}\label{4.11.1}\begin{array}{l}\xi =\xi _1+2^3\xi _2+2^{p-7}\xi _4\\
\ \\
0<\xi _1<2^3,\ 0\le \xi _2<2^{p-10}, 0\le \xi _4<2^6 \end{array}\end{equation}
We have 
\begin{equation}\label{4.11.2}\lambda (k_{t-1},k_t,p_t,p_{[t+1,r]}, \xi /2^p)\le \lambda _0\end{equation}
 if there is $\ell $ such that $k_{j-1}<\ell <\ell +M<k_j$ and $x_{[k_{j-1},\ell +M]} <p_j$ and for all $1\le u\le M$ and $0\le v\le M$, and $n=x_{[\ell,k]}$ the following {\em{does not hold}}:
\begin{equation}\label{4.11.3}\begin{array}{l}a_{n+p+u,\ell +v}=a_{n+p,\ell},\ 0\le u,v\le M\\
\ \\
0\le \xi _4<2^5\mbox{ or }a_{n+u,\ell +v}=a_{n+1,\ell },\ 0\le u,v\le M\end{array}\end{equation}
\end{theorem}

\noindent{\em{Proof.}} This is the analysis in \ref{4.9}.  $\lambda (k_{t-1},k_t,p_t,p_{[t+1,r]}, \xi )$ is a convex sum of roots of unity where the number of roots of unity is 
$$\left(\begin{matrix}p_t\\ k_t-k_{t-1}-1\end{matrix}\right) $$
and thus bounded in terms of $M'$ and the coefficient of each root of unity is the inverse of this. Since $p_t$ is bounded there are only finitely many values for $(x_{[k_{t-1}+1,k]},\cdots x_{[\ell,k]}\cdots x_{[k_t,k]})$ and again the number of values is bounded in terms of $M'$. But given the condition $M'\ge M^2+p_t-(k_t-k_{t-1})$, there will be 
at least one $(x_{[k_{t-1}+1,k]},\cdots x_{[\ell,k]}\cdots x_{[k_t,k]})$  with an $\ell\in(k_{t-1},k_t)$ such that $x_{[\ell +v,k]}>x_{[\ell +v+1,k]}+M$ for $0\le v\le M$. So the analysis of \ref{4.9} can be applied with $n=x_{[\ell +v,k]}$ for each $0\le v<M$.
\ep

So then Conjecture \ref{4.7} is implied by Conjecture \ref{4.12}.

\begin{conjecture}\label{4.12} For $1<\mu <\infty $ and suitable constants  $M>0$, and $c>0$, and for $k$ sufficiently large,  and $p\le k$: apart from a set of  $[x_1,\cdots x_k]$ of $\mathbb G_{1-\mu ^{-1}}$-measure $<e^{-ck}$ the set of $\ell \le k$ such that (\ref{4.12.1}) holds is of density $\le (1-c)k$, where $n=x_{[\ell ,k]}$, 
\begin{equation}\label{4.12.1}a_{n+p+u,\ell +v}=a_{n+p,\ell} \mbox{  for }1\le u\le M,\ 0\le v\le M\end{equation}
\end{conjecture}

This can be simplified  as follows - where $c$ does not transfer form \ref{4.12} to \ref{4.13}.

\begin{conjecture}\label{4.13}  For $1<\mu <\infty $ and  a suitable constant  $c>0$, and $k$ sufficiently large and  $p\le k$ apart from a set of  $[x_1,\cdots x_k]$ of $\mathbb G_{1-\mu ^{-1}}$-measure $<e^{-ck}$ the set of $\ell \le k$ such that $a_{x_{[\ell,k]}+p,\ell }=0$ has density $\le (1-c)k$ and similarly for the  set of $\ell \le k$ such that $a_{x_{[\ell,k]}+p,\ell }=1$.\end{conjecture}

\subsection{Comments on Conjecture \ref{4.13}}\label{4.14}
 The parallel with Section 7 of \cite{Tao} is striking. We have an infinite array with entries labelled $(j,\ell )$  for $j\ge 0$ and $\ell \ge 1$. and we also have random paths on this infinite array.  It is convenient to label the axes so that the horizontal axis labelled $j$ runs from left to right starting from $j=0$and the vertical axis labelled $\ell $ points downwards,  starting from $\ell =1$. This is because one usually starts writing from the top of the page and in roman script from left to right. The $\ell$'th row of the array is the coefficients $a_{\ell ,j}\in \{ 0,1\} $ where
$$-3^{-\ell }=\sum _{j=0}^\infty a_{j,\ell }2^\ell .$$
The paths $(\ell, x_{[\ell ,k]})$ in \ref{4.12} run upwards and to the right -- which actually is also the case in \cite{Tao} for paths $(j,b_{[1,j]})$. But the analogues of the  ``black triangles'' in \cite{Tao}, ``constant triangles'' of $0$'s or of $1$'s, which will appear in Part 2 (hopefully), will be differently oriented because of the different axes direction.   The array in \cite{Tao} is somewhat more closely associated with triadic expansions of $2^{-n}$, while here we are working explicitly  with dyadic expansions of $-3^{-k}$. Actually as we shall see in Part 2, dyadic expansions of $-3^{-k}$ and triadic expansions of $2^{-n}$ are closely related,  and it is helpful to consider them together, with yet another axis change and an affine transformation of constant triangles.

 \end{document}